\documentclass[11pt]{article}
\usepackage{amssymb,amsmath, amsfonts, amsthm, amssymb}
\usepackage{bbding}
\usepackage[dvips]{graphicx}
\graphicspath{{./figuras/}}
\usepackage{subfigure}

\usepackage[total={6.5in,9.00in}, top=1.2in, left=0.9in, includefoot]{geometry}

\newcommand{\norm}[1]{\left\Vert#1\right\Vert}
\newcommand{\seq}[1]{\left<#1\right>}
\newcommand{\abs}[1]{\left\vert#1\right\vert}
\newcommand{\RE}{{\rm I}\!{\rm R}} 
\newcommand{\N}{{\rm I}\!{\rm N}} 

\newcommand{\et}{\mathcal{E}_{\{R_{\alpha}\}}^{\text{\rm{tot}}}}

 \newtheorem{thm}{Theorem}[section]
 
 \newtheorem{lem}[thm]{Lemma}
 
 \newtheorem{defn}[thm]{Definition}
 \newtheorem{rem}[thm]{Remark}
 \newtheorem{ej}[thm]{Example}
\title{On the existence of global saturation for spectral regularization methods with optimal
qualification\thanks{This work was supported in part by Consejo
Nacional de Investigaciones Cient\'{\i}ficas y T\'{e}cnicas,
CONICET, through PIP 2010-2012 Nro. 0219, by Universidad Nacional
del Litoral, U.N.L., through project CAI+D 2009-PI-62-315, by
Agencia Nacional de Promoci\'{o}n Cient\'{\i}fica y Tecnol\'{o}gica,
ANPCyT, through project PICT-2008-1301 and by the Air Force Office
of Scientific Research, AFOSR, through Grant FA9550-10-1-0018.}}

\author{
Gisela L. Mazzieri{\footnote{\textrm{Instituto de Matem\'{a}tica
Aplicada del Litoral, IMAL, CONICET-UNL, G\"{u}emes 3450, S3000GLN,
Santa Fe, Argentina.}}}
\thanks{Departamento de Matem\'{a}tica, Facultad de
Bioqu\'{\i}mica y Ciencias Biol\'{o}gicas, Universidad Nacional del
Litoral, Santa Fe, Argentina({\tt gmazzieri@hotmail.com}).}
\and
Ruben D. Spies$^{\text{\Envelope\,\,}
\dagger}$\thanks{Departamento de Matem\'{a}tica, Facultad de
Ingenier\'{\i}a Qu\'{\i}mica, Universidad Nacional del Litoral, Santa Fe,
Argentina (\Envelope\,: {\tt rspies@santafe-conicet.gov.ar}).}
\and
Karina G. Temperini$^{\dagger}$\thanks{Departamento de Matem\'{a}tica,
Facultad de Humanidades y Ciencias, Universidad Nacional del
Litoral, Santa Fe, Argentina ({\tt
ktemperini@santafe-conicet.gov.ar}).} }
\date{\empty}
\begin{document}
\maketitle

\begin{abstract}
%

A family of real functions $\{g_\alpha\}$ defining a spectral
regularization method with optimal qualification is considered.
Sufficient condition on the family and on the optimal
qualification guaranteeing the existence of saturation are
established. Appropriate characterizations of both the saturation
function and the saturation set are found and some examples are
provided.
\end{abstract}



\pagestyle{myheadings} \thispagestyle{plain} \markboth{GISELA
MAZZIERI, RUBEN SPIES and KARINA TEMPERINI }{EXISTENCE OF SATURATION FOR
SPECTRAL REGULARIZATION METHODS}

\section{Introduction}\label{intro}
Let $X, Y$ be infinite dimensional Hilbert spaces and
$T:X\rightarrow Y$ a bounded linear operator with non-closed range
$\mathcal{R}(T)$. It is well known that under these conditions
$T^\dag$, the Moore-Penrose generalized inverse of $T$, is unbounded
(\cite{bookEHN}) and therefore the linear operator equation
\begin{equation}\label{eq:prob}
    Tx=y
\end{equation}
is ill-posed. The
Moore-Penrose generalized inverse can be used to define the least
squares solutions of (\ref{eq:prob}). In fact equation (\ref{eq:prob}) has a
least squares solution if and only if $y \in
\mathcal{D}(T^\dag)\doteq \mathcal{R}(T)\oplus
\mathcal{R}(T)^\perp$ and in that case, $x^\dag\doteq T^\dag y$ is the
least squares solution of minimum norm and the set of all
least-squares solutions of (\ref{eq:prob}) is given by $x^\dag +
\mathcal{N}(T)$. Since $T^\dag$ is unbounded, $x^\dag$ does not
depend continuously on the data $y$. Therefore, if instead of the
exact data $y$, a noisy observation $y^\delta$ is available,
$y^\delta=Tx+\delta \xi$, where the noise $\xi$ is assumed to be
bounded, $\norm{\xi}\leq 1$, then it is possible that $T^\dag
y^\delta$ does not even exist and if it does, it will not
necessarily be a good approximation of $x^\dag$ (\cite{ref:Seidman-80}, \cite{Spies-Temperini-2006}). This instability
becomes evident when trying to approximate $x^\dag$ by traditional
numerical methods and procedures.

Ill-posed problems must be first regularized if one wants to
successfully attack the task of numerically approximating their
solutions. Regularizing an ill-posed problem such as
(\ref{eq:prob}) essentially means approximating the operator
$T^\dag$ by a parametric family of bounded operators
$\{R_\alpha\}$, where $\alpha$ is the so called ``regularization
parameter''. If $y \in \mathcal{D}(T^\dag)$, then the best
approximate solution $x^\dag$ of (\ref{eq:prob}) can be written as
$x^\dag=\int_0^{\norm{T}^2+} \frac{1}{\lambda} \, dE_\lambda
T^\ast y$ where $\{E_\lambda\}$ is the spectral family associated
to the operator $T^\ast T$ (see \cite{bookEHN}). This is mainly
why many regularization methods are based on spectral theory and
consist of defining $R_\alpha\doteq \int_0^{\norm{T}^2+}
g_\alpha(\lambda) \, dE_\lambda T^\ast$ where $\{g_\alpha\}$ is a
family of functions appropriately chosen such that for every
$\lambda \in (0,\norm{T}^2]$ there holds $\underset{\alpha
\rightarrow 0^+}{\lim}g_\alpha(\lambda)=\frac{1}{\lambda}$. It is
important to emphasize however that no mathematical trick can make
stable a problem that is intrinsically unstable. Whatever the
case, there is always loss of information. All a regularization
method can do is to recover the largest possible amount of
information about the solution of the problem, maintaining
stability. It is often said that the art of applying
regularization methods consist always in maintaining an adequate
balance between accuracy and stability. Usually accuracy can be
improved with increasing assumptions (or information) on the
regularity of the exact solution. In 1994, however, Neubauer
(\cite{Neubauer94}) showed that certain spectral regularization
methods ``{\it saturate}'', that is, they become unable to
continue extracting additional information about the exact
solution even upon increasing regularity assumptions on it. In his
article, Neubauer introduced for the first time the idea of the
concept of ``{\it saturation}'' of regularization methods.
Saturation is however a rather subtle and complex issue in the
study of regularization methods for inverse ill-posed problems and
the concept has, for many years, escaped rigorous formalization in
a general context. Neubauer's idea referred to the best order of
convergence that a method can achieve independently of the
smoothness assumptions on the exact solution and on the selection
of the parameter choice rule. In 1997, Neubauer
(\cite{ref:Neubauer-97}) showed that this saturation phenomenon
occurs for instance in the classical Tikhonov-Phillips method.
Later on, in 2004, Math\'{e} (\cite{Mathe2004}) proposed a general
definition of the concept of saturation for spectral
regularization methods. However, the concept of saturation defined
by Math\'{e} is not applicable to general regularization methods and
it is not fully compatible with the original idea of saturation
proposed by Neubauer in \cite{Neubauer94}. In particular, for
instance, the definition of saturation given in \cite{Mathe2004}
does not imply uniqueness and therefore, neither a best global
order of convergence. More recently, in 2011, Herdman, Spies and
Temperini (see \cite{ref:Herdman-Spies-Temperini-2011}) developed
a general theory of global saturation for arbitrary regularization
methods, formalizing the original and intuitive idea first
outlined by Neubauer in 1994 (\cite{Neubauer94}).

Related in a dual way to the concept of saturation is the concept
of qualification of a spectral regularization method, introduced
by Math\'{e} and Pereverzev in 2003
(\cite{ref:Mathe-Pereverzev-2003}). This concept is strongly
related to the optimal order of convergence of the regularization
error, under certain ``a-priori'' assumptions on the exact solution.
In 2009 Herdman, Spies and Temperini
(\cite{ref:Herdman-Spies-Temperini-2009}) generalized the concept
of qualification and introduced three hierarchical levels of it:
weak, strong and optimal qualification. There, it was shown that the weak
qualification generalizes the definition introduced in
\cite{ref:Mathe-Pereverzev-2003}.

In this work, some light on the existence of saturation for
spectral regularization methods with optimal qualification is
shed. In particular, sufficient conditions on the family of real
functions $\{g_\alpha\}$ defining the method and on the optimal
qualification $\rho$, which guarantee the existence of saturation,
are established. Moreover, in those cases, appropriate
characterizations of both the saturation function and the
saturation set are provided.

\section{Preliminaries} \label{sec:2}
In this section we shall recall some basic concepts on global
saturation of regularization methods for inverse ill-posed
problems theory (for more details see
\cite{ref:Herdman-Spies-Temperini-2011}). In the sequel,
$T:X\rightarrow Y$ will be a bounded linear operator with
non-closed range between two Hilbert spaces $X$ and $Y$. Without
loss of generality we will assume that the operator $T$ is
invertible (in the context of inverse problems it is customary to
work with the Moore-Penrose generalized inverse of $T$ since one
seeks least squares solutions of the problems; therefore the lack
of injectivity of $T$ is never a relevant issue). Also, for
simplicity of notation and unless otherwise specified, we shall
assume that all subsets of the Hilbert space $X$ under
consideration are not empty and they do not contain $x=0$.

Let $M\subset X$. We shall say that a function $\psi:X \times \RE \to
\RE$ belongs to the class $\mathcal{F}_{M}$ if there exists
$a=a(\psi)>0$ such that $\psi$ is defined in $M \times (0,a)$, with
values in $(0,\infty)$ and it satisfies the following conditions:
\begin{enumerate}
    \item $\underset{\delta\rightarrow 0^+}{\lim}\psi(x,\delta)=0$
    for all
$x \in
    M$, and
    \item  $\psi$ is continuous and non-decreasing as a function of $\delta$
    in $(0,a)$ for each fixed $x\in M$.
\end{enumerate}
\bigskip One may think of $\mathcal{F}_{M}$ as the collection of all possible $\delta$-``{\it orders of convergence}'' on
the set $M$.
\medskip
\begin{defn}
Let $\{R_{\alpha}\}_{\alpha\in (0,\alpha_0)}$ be a family of
regularization operators for the problem $Tx=y$. The ``total error
of $\{R_{\alpha}\}_{\alpha\in (0,\alpha_0)}$ at $x \in X$ for a
noise level $\delta$'' is defined as
$$\et(x,\delta)\doteq
\underset{\alpha \in (0,\alpha_0)}{\inf}\;\underset{y^\delta \in
\overline{B_\delta(Tx)}}{\sup} \norm{R_\alpha y^\delta-x},$$ where
$\overline{B_\delta(Tx)}\doteq\{y \in Y: \norm{Tx-y}\leq \delta\}$.
\end{defn}
\smallskip

Note that $\et$ is the error in the sense of the largest possible
discrepancy that can be obtained for an observation of $y$ within
noise level $\delta$, with an appropriate choice of the
regularization parameter $\alpha$.
\medskip

\begin{defn}\label{def:relac-fm}
Let $M\subset X$ and $\psi,\tilde{\psi} \in \mathcal{F}_{M}$. We say
that ``$\psi$ precedes $\tilde{\psi}$ on $M$'', and we denote it
with $\psi \overset{M}{\preceq} \tilde{\psi}$, if there exist a
constant $r>0$ and a function $p:M\rightarrow (0,\infty)$ such that
$\psi(x,\delta) \leq p(x)\, \tilde{\psi}(x,\delta)$ for all $x \in
M$ and for every $\delta \in (0,r)$.
\end{defn}
\medskip

\begin{defn}\label{def-sat}
Let $\{R_{\alpha}\}_{\alpha \in (0,\alpha_0)}$ be a family of
regularization operators for the problem $Tx=y$, $M\subset X$ and
$\psi \in \mathcal{F}_M$. We say that $\psi$ is an ``upper bound of
convergence for the total error of $\{R_\alpha\}_{\alpha \in
(0,\alpha_0)}$ on $M$'' if $\et\overset{M}{\preceq}\psi$.
\end{defn}

\smallskip

With $\mathcal{U}_M(\et)$ we shall denote the set of all functions
$\psi \in \mathcal{F}_M$ that are upper bounds of convergence for
the total error of $\{R_\alpha\}_{\alpha \in (0,\alpha_0)}$ on $M$.

The following two definitions formalize certain comparisons of
bounds of convergence on different sets of $X$, which will be needed
later to introduce the concept of global saturation.

\begin{defn}\label{def:relac-mn}
Let $M,\tilde{M}\subset X$, $\psi \in \mathcal{F}_M$ and
$\tilde{\psi}\in \mathcal{F}_{\tilde{M}}$.

i) We say that ``$\psi$ on $M$ precedes $\tilde{\psi}$ on
$\tilde{M}$'', and we denote it with $\psi \overset{M,
\tilde{M}}{\preceq}\tilde{\psi}$, if there exist a constant $d>0$
and a function $k:M\times \tilde{M}\rightarrow (0,\infty)$ such that
$\psi(x,\delta)\leq k(x,\tilde{x})\, \tilde{\psi}(\tilde{x},\delta)$
for every $x \in M$, $\tilde{x}\in \tilde{M}$ and $\delta \in
(0,d).$

ii) We say that ``$\psi$ on $M$ is equivalent to $\tilde{\psi}$ on
$\tilde{M}$'', and we denote it with $\psi \overset{M,
\tilde{M}}{\approx}\tilde{\psi}$, if $\,\psi \overset{M,
\tilde{M}}{\preceq}\tilde{\psi}\;$ and $\;\tilde{\psi}
\overset{\tilde{M}, M}{\preceq}{\psi}$.
\end{defn}

\medskip
\begin{defn}
Let $M\subset X$ and $\psi \in \mathcal{F}_M$. We say that ``$\psi$
is invariant over $M$'' if $\psi \overset{M,M}{\approx}\psi$.
\end{defn}
\medskip

Next we recall the concept of global saturation introduced in
\cite{ref:Herdman-Spies-Temperini-2011}.

\begin{defn}\label{def:satur}
Let $M_S \subset X$ and $\psi_S \in \mathcal{U}_{M_S}(\et)$. It is
said that $\psi_S$ is a ``global saturation function of
$\{R_\alpha\}$ over $M_S$'' if $\psi_S$ satisfies the following
three conditions:

S1. For every $x^\ast \in X$, $x^\ast\neq 0$, $x \in M_S$,
$\underset{\delta
\rightarrow0^+}{\limsup}\,\frac{\et(x^\ast,\delta)}{\psi_S(x,\delta)}>0.$

S2. $\psi_S$ is invariant over $M_S$.

S3. There is no upper bound of convergence for the total error of
$\{R_\alpha\}$ that is a proper extension of $\psi_S$ (in the
variable $x$) and satisfies S1 and S2, that is, there exist no
$\tilde{M}\supsetneqq M_S$ and $\tilde{\psi} \in \mathcal{U}_{
\tilde{M}}(\et)$ such that $\tilde{\psi}$ satisfies S1 and S2 with
$M_S$ replaced by $\tilde{M}$ and $\psi_S$ replaced by
$\tilde{\psi}$.
\end{defn}
\smallskip

The function $\psi_S$ and the set $M_S$ are refer to as the
saturation function and the saturation set, respectively.

This conception of global saturation essentially establishes that in
no point $x^\ast \in X$, $x^\ast \neq 0$, can exist an upper bound
of convergence for the total error of the regularization method that
is ``strictly better'' than the saturation function $\psi_S$ at any
point of the saturation set $M_S$.

\medskip
Let $\{E_\lambda\}_{\lambda \in \RE}$ be the spectral family
associated to the linear selfadjoint operator $T^\ast T$ and
$\{g_\alpha\}_{\alpha \in (0, \alpha_0)}$ a parametric family of
functions  $g_\alpha:[0,\norm{T}^2]\rightarrow \RE$, $\alpha \in
(0,\alpha_0)$, and consider the following standing hypotheses:

\textit{H1}.\; For every $\alpha \in (0,\alpha_0)$ the function
$g_\alpha$ is piecewise continuous on $[0,\norm{T}^2]$.

\textit{H2}.\; There exists a constant $C>0$ (independent of
$\alpha$) such that $\abs{\lambda g_\alpha (\lambda)}\leq C$ for
every $\lambda \in [0,\norm{T}^2]$.

\textit{H3}.\; For every  $\lambda \in (0,\norm{T}^2]$, there exists
$\underset{\alpha \rightarrow 0^+}{\lim}
g_\alpha(\lambda)=\frac{1}{\lambda}$.

\smallskip

If $\{g_\alpha\}_{\alpha \in (0, \alpha_0)}$ satisfies hypotheses
\textit{H1-H3}, then (see \cite{bookEHN}, Theorem 4.1) the
collection of operators $\{R_\alpha\}_{\alpha \in (0, \alpha_0)}$,
where
\begin{equation*}
    R_\alpha\doteq \int_0^{\norm{T}^2+} g_\alpha(\lambda)\, dE_\lambda \,T^\ast
    =g_\alpha(T^\ast T)T^\ast,
\end{equation*}
is a family of regularization operators for $T^\dag$. In this case
we say that $\{R_\alpha\}_{\alpha \in (0, \alpha_0)}$ is a ``family
of spectral regularization operators'' (FSRO) for $Tx=y$ and
$\{g_\alpha\}_{\alpha \in (0, \alpha_0)}$ is a ``spectral
regularization method'' (SRM).

The following definitions will be needed both to recall the concept
of qualification as introduced in
\cite{ref:Herdman-Spies-Temperini-2009}, as well as in the rest of
the article.

We denote with $\mathcal{O}$ the set of all non-decreasing functions
$\rho:\RE^+ \to \RE^+$ such that $\underset{\alpha \rightarrow
0^+}{\lim}\rho(\alpha)=0$ and with $\mathcal{S}$ the set of all
continuous functions $s:\RE^+_0 \to \RE^+_0$ satisfying $s(0)=0$ and
such that $s(\lambda)>0$ for every $\lambda>0.$ Note that if $s\in
\mathcal{S}$ is non-decreasing, then $s$ is an \textit{index
function} in the sense of Math\'{e}-Pereverzev
(\cite{ref:Mathe-Pereverzev-2003}).

\begin{defn}
Let $\rho,\tilde{\rho}\in \mathcal{O}$. We say that ``$\rho$
precedes $\tilde{\rho}$ at the origin'' and we denote it with $\rho
\preceq \tilde{\rho}$, if there exist positive constants $c$ and
$\varepsilon$ such that $\rho(\alpha)\leq c\,\tilde{\rho}(\alpha)$
for every $\alpha \in (0,\varepsilon)$.
\end{defn}

\begin{defn} \label{def:order}
Let $\rho,\tilde{\rho}\in \mathcal{O}$. We say that ``$\rho$ and
$\tilde{\rho}$ are equivalent at the origin'' and we denote it with
$\rho \approx\tilde{\rho}$, if they precede each other at the
origin, that is, if there exist constants $\varepsilon>0,\, c_1,\,
c_2$, $0<c_1<c_2<\infty$ such that $c_1\, \rho(\alpha)\leq
\tilde{\rho}(\alpha) \leq c_2\, \rho(\alpha)$ for every $\alpha \in
(0,\varepsilon)$.
\end{defn}

Clearly ``$\approx$'' is an equivalence relation and it introduces
in $\mathcal{O}$ a partial ordering. Analogous definitions and
notation will be used for $s,\tilde{s} \in \mathcal{S}$.

\begin{defn}
Let $\{g_\alpha\}_{\alpha \in (0,\alpha_0)}$ be a SRM,
$r_{\alpha}(\lambda)\doteq 1-\lambda g_{\alpha}(\lambda)$, $\rho\in
\mathcal{O}$ and $s \in \mathcal{S}.$

\textbf{i)} We say that $(s,\rho)$ is a ``weak source-order pair for
$\{g_\alpha\}$'' if it satisfies
\begin{equation}\label{eq:O}
\frac{s(\lambda)\abs{r_\alpha(\lambda)}}{\rho(\alpha)}=O(1)\quad
\textrm{for}\;\; \alpha \rightarrow 0^+, \; \forall\; \lambda>0.
\end{equation}

\textbf{ii)} We say that $(s, \rho)$ is a ``strong source-order pair
for $\{g_\alpha\}$'' if it is a weak source-order pair and there is
no $\lambda>0$ for which ``$O(1)$'' in (\ref{eq:O}) can be replaced
by ``$o(1)$''. That is, if $(s,\rho)$ is a weak source-order pair
for $\{g_\alpha\}$ and
\begin{equation}\label{eq:no-o}
\underset{\alpha \rightarrow
0^+}{\limsup}\,\frac{s(\lambda)\abs{r_\alpha(\lambda)}}{\rho(\alpha)}>0\quad
\forall\; \lambda>0.
\end{equation}

\textbf{iii)} We say that $(\rho, s)$ is an ``order-source pair for
$\{g_\alpha\}$'' if there exist a constant $\gamma>0$ and a function
$h:(0,\alpha_0)\rightarrow \RE^+$ with $\underset{\alpha
\rightarrow0^+}{\lim}h(\alpha)=0$, such that
\begin{equation}\label{eq:4.50}
\frac{s(\lambda)\abs{r_\alpha(\lambda)}}{\rho(\alpha)}\geq \gamma
\quad \forall\; \lambda \in [h(\alpha),+\infty).
\end{equation}
\end{defn}

In the context of the previous definitions we refer to the function
$\rho$ as the ``order of convergence'' and to $s$ as the ``source
function''.

We are now ready to define the concept of qualification in its three
different levels as it was introduced in
\cite{ref:Herdman-Spies-Temperini-2009}.

\begin{defn}\label{def:calif-3} Let $\{g_\alpha\}$ be a SRM.

\textbf{i)} We say that $\rho$ is ``weak qualification of
$\{g_\alpha\}$'' if there exists a function $s$ such that $(s,\rho)$
is a weak source-order pair for $\{g_\alpha\}$.

\textbf{ii)} We say that $\rho$ is ``strong qualification of
$\{g_\alpha\}$'' if there exists a function $s$ such that $(s,\rho)$
is a strong source-order pair for $\{g_\alpha\}$.

\textbf{iii)} We say that $\rho$ is ``optimal qualification of
$\{g_\alpha\}$'' if there exists a function $s$ such that $(s,\rho)$
is a strong source-order pair for $\{g_\alpha\}$ and $(\rho, s)$ is
an order-source pair for $\{g_\alpha\}$.
\end{defn}

Note that since condition (\ref{eq:4.50}) implies condition
(\ref{eq:no-o}), in the definition of optimal qualification above
the requirement that $(s,\rho)$ be strong source-order pair can be
replaced by the one that $(s,\rho)$ be a weak source-order pair.

Now given the SRM $\{g_\alpha\}$ and $\rho\in \mathcal{O}$, we
define
\begin{equation}\label{s-rho}
s_\rho(\lambda)\doteq  \underset{\alpha \to 0^+}{\liminf}
\frac{\rho(\alpha)}{\abs{r_\alpha(\lambda)}}\quad \textrm{for} \quad
\lambda \geq 0.
\end{equation}
Note that $s_\rho(0)=0$ and if $s_\rho$ is continuous, $s_\rho \in
\mathcal{S}.$

The next theorem provides necessary and sufficient condition, in
terms of $s_\rho$, for an order of convergence $\rho \in
\mathcal{O}$ to be optimal qualification.

\begin{thm}\label{teo:cond-calopt}
(\cite{ref:Herdman-Spies-Temperini-2009}) Let $\{g_\alpha\}$ be a
SRM and $\rho \in \mathcal{O}$ such that $s_\rho \in \mathcal{S}$.
Then $\rho$ is optimal qualification of $\{g_\alpha\}$ if and only
if $s_\rho$ verifies \text{\rm (\ref{eq:4.50})} and
\begin{equation}\label{eq:cond-calif}
0<s_\rho(\lambda)<+\infty\quad  \textrm{for every}\quad \lambda
>0.
\end{equation}
\end{thm}

The next theorem shows the uniqueness of the source function.
\begin{thm}\label{teo:unicaff} (\cite{ref:Herdman-Spies-Temperini-2009})
If $\rho$ is optimal qualification of $\{g_\alpha\}$ then there
exists only one function $s$ (in the sense of the equivalence
classes induced by Definition \ref{def:order}) such that $(s,\rho)$
is a strong source-order pair and $(\rho,s)$ is an order-source pair
for $\{g_\alpha\}$. Moreover if $s_\rho \in \mathcal{S}$, then
$s_\rho$ is such a unique function.
\end{thm}

The following converse result, where regularity properties of the
exact solution are derived in terms of the rate of convergence of
the regularization error, will be needed later. This result states
that if the regularization error has order of convergence
$\rho(\alpha)$ and $(\rho,s)$ is an order-source pair, then the
exact solution belongs to the source set given by the range of the
operator $s(T^\ast T)$.

\begin{thm}\label{teo:gen-411}
(\cite{ref:Herdman-Spies-Temperini-2009}) Let $\{g_\alpha\}$ be a
SRM and $R_\alpha=g_\alpha(T^\ast T)T^\ast$. If
$\norm{(R_\alpha-T^\dag)y}=O(\rho(\alpha))$ for $\alpha\rightarrow
0^+$ and $(\rho,s)$ is an order-source pair for $\{g_\alpha\}$,
then $T^\dag y\in \mathcal{R}(s(T^\ast T)).$
\end{thm}

\section{Saturation of spectral regularization methods with optimal qualification} \label{sec:3}
The purpose of this section is to shed some light on the saturation of SRM with
optimal qualification. More precisely, we will establish sufficient
conditions on the family of functions $\{g_\alpha\}$ and on the
optimal qualification $\rho$ guaranteeing the existence of
saturation. Moreover, for those methods we will provide appropriate
characterizations of both the saturation function and the saturation
set. Then, let $\{g_\alpha\}_{\alpha \in (0,\alpha_0)}$ be a SRM and consider the following hypothesis:

\textit{H4}.\; Exists $k>0$ such that $G_\alpha \doteq
\norm{g_\alpha(\cdot)}_\infty\leq \frac{k}{\sqrt{\alpha}} \; \forall \,\alpha \in (0,\alpha_0).$

\begin{lem}\label{lem:cotasup}
Let $\{g_\alpha\}_{\alpha\in (0,\alpha_0)}$ be a SRM satisfying hypothesis \textit{H4}, $R_\alpha=g_\alpha(T^\ast T)T^\ast$ and
$(s,\rho)$ a weak source-order pair for $\{g_\alpha\}_{\alpha\in
(0,\alpha_0)}$ where $\rho$ is continuous. Define $X^s\doteq \mathcal{R}(s(T^\ast
T))\setminus \{0\}$, $\Theta(t)\doteq \sqrt{t}\,\rho(t)$ for $t>0$, $\psi(x,\delta)\doteq \rho \circ \Theta^{-1}(\delta)$ for $x
\in X^s$ and $\delta \in (0, \Theta(\alpha_0))$. Then $\psi \in \mathcal{F}_{X^s}$ and, moreover, $\psi$ is an upper bound of
convergence for the total error of $\{R_\alpha\}_{\alpha\in
(0,\alpha_0)}$ on $X^s$, that is, $\psi \in \mathcal{U}_{X^s}(\et)$.
\end{lem}

\smallskip

\begin{proof}
Since $\rho$ is continuous and non-decreasing and $\rho(0^+)=0$ it follows that $\Theta(t)$ is continuous and strictly increasing on $(0, +\infty)$ with $\Theta(0^+)=0.$ Therefore $\Theta^{-1}$ exists and has the same properties. It then follows that $\psi$ is continuous and non-decreasing as a function of $\delta$ in $(0, \Theta(\alpha_0))$ for each fixed $x\in X^s$, and $\psi(x,0^+)=0$ for all $x\in X^s$. Hence $\psi \in \mathcal{F}_{X^s}.$

On the other hand, since $(s,\rho)$ is a weak source-order pair for
$\{g_\alpha\}_{\alpha\in (0,\alpha_0)}$, there exist positive
constants $c$ and $\hat{\alpha}$ such that $s(\lambda)\abs{r_\alpha(\lambda)}\leq c\rho(\alpha)$ for every $\alpha \in (0,\hat{\alpha}), \lambda \in (0,\norm{T}^2].$ Moreover, from hypothesis \textit{H2} and the fact that $\rho$ is non-decreasing it follows that the previous inequality holds (perhaps with a different positive constant $c$) for every $\alpha \in (0,\alpha_0),$ that is,
\begin{equation}\label{eq:weak}
s(\lambda)\abs{r_\alpha(\lambda)}\leq c\rho(\alpha)\quad \forall \;
\alpha \in (0,\alpha_0), \forall\; \lambda \in (0,\norm{T}^2].
\end{equation}

Now, for every $p\geq 0$ we define the source sets $X^{s,p}\doteq
\{x \in X:x=s(T^\ast T)\zeta$, $\norm{\zeta}\leq p \}$. Then for each $x \in
X^{s}$ there exists $p_x > 1$ such that $x \in X^{s,p_x}$. On the other hand, since
$\Theta$ is continuous and strictly increasing in $(0,\alpha_0)$,
there exists a unique $\tilde{\alpha}_x \in (0,\alpha_0)$ such that
$x\in X^{s,p_x}$ and $\Theta(\tilde{\alpha}_x)=\frac{\delta}{p_x}.$
Therefore,
\begin{eqnarray}\label{eq:Trho}\nonumber
  \et(x,\delta) &=& \underset{\alpha \in
(0,\alpha_0)}{\inf}\;\underset{y^\delta \in
  \overline{B_\delta(Tx)}}{\sup}\norm{R_\alpha\,
y^\delta-x} \\ \nonumber
   &\leq & \underset{y^\delta \in
   \overline{B_\delta(Tx)}}{\sup}\norm{R_{\tilde{\alpha}_x}\,
y^\delta - x}.
\end{eqnarray}
Now, since $y^\delta=Tx+\delta \xi$, $\norm{\xi}\leq 1$ and
$x=s(T^\ast T)\zeta$ with $\norm{\zeta}\leq p_x$, it follows
immediately that
\begin{eqnarray}\label{eq:sup}\nonumber
\norm{R_{\tilde{\alpha}_x} \,y^\delta-x} &\leq&
\norm{(g_{\tilde{\alpha}_x}(T^\ast T)T^\ast T-I)s(T^\ast
T)\zeta}+\delta\norm{g_{\tilde{\alpha}_x}(T^\ast T)T^\ast\xi}\\
 &\leq& p_x \underset{\lambda \in
 (0,\norm{T}^2]}{\sup}\{s(\lambda)\abs{r_{\tilde{\alpha}_x}(\lambda)}\}+\delta\underset{\lambda \in
 (0,\norm{T}^2]}{\sup}\{\sqrt{\lambda}\abs{g_{\tilde{\alpha}_x}(\lambda)}\},
\end{eqnarray}
where the last inequality follows from properties of functions of a
selfadjoint operator, (more precisely, for any piecewise continuous function $f$ there holds
$\norm{f(T^\ast T)}\leq \underset{\lambda}{\sup}\abs{f(\lambda)}$ and
$\norm{f(T^\ast T)T^\ast}\leq
\underset{\lambda}{\sup}\{\sqrt{\lambda}\abs{f(\lambda)}\}$, see \cite{bookEHN}, p. 45). Using
(\ref{eq:weak}) and hypothesis \textit{H4} in (\ref{eq:sup}) it follows that
\begin{equation}\label{eq:des}
\norm{R_{\tilde{\alpha}_x} \,y^\delta-x}\leq
p_x\,c\,\rho(\tilde{\alpha}_x)+\delta \frac{k}{\sqrt{\tilde{\alpha}_x}}
\norm{T}.
\end{equation}
Since
$\Theta(\tilde{\alpha}_x)=\sqrt{\tilde{\alpha}_x}\,\rho(\tilde{\alpha}_x)=\frac{\delta}{p_x}$,
it follows that
$\frac{\delta}{\sqrt{\tilde{\alpha}_x}}=p_x\,\rho(\tilde{\alpha}_x).$
Hence by virtue of (\ref{eq:des}) one has that
\begin{eqnarray}\label{eq:s}\nonumber
 \norm{R_{\tilde{\alpha}_x} \,y^\delta-x}&\leq& p_x(c
+k\norm{T})\rho(\tilde{\alpha}_x)\\ \nonumber
&=&p_x(c+k\norm{T})\rho\left(\Theta^{-1}\left(\frac \delta
p_x\right)\right)\\
&\leq & p_x(c+k\norm{T})\rho(\Theta^{-1}(\delta)),
\end{eqnarray}
where the last inequality follows from the fact that $p_x> 1$ and
both $\rho$ and $\Theta^{-1}$ are non-decreasing functions. From
(\ref{eq:Trho}) and (\ref{eq:s}) it follows that for every $\delta
\in (0,\Theta(\alpha_0))$,
\begin{equation*}
    \et(x,\delta)\leq p_x(c
+k\norm{T}) \rho(\Theta^{-1}(\delta)) =p_x(c
+k\norm{T})\,\psi(x,\delta).
\end{equation*}
This proves that $\psi \in \mathcal{U}_{X^s}(\et)$.
\end{proof}

\begin{defn}\label{def:local-type} (\cite{ref:Iaffei-2005})
Let $k$ be a positive constant and
$\rho:(0,k]\rightarrow(0,+\infty)$ a continuous non-decreasing
function such that $\underset{t\rightarrow 0^+}{\lim}\rho(t)=0.$
We say that $\rho$ is of local upper type $\beta$ ($\beta \geq 0$)
if there exists a positive constant $d$ such that $\displaystyle
\rho(t)\leq ds^{-\beta}\rho(s\,t)$ for every $s \in (0,1]$, $t \in
(0,k]$.
\end{defn}

\smallskip
\begin{thm}\label{teo:sat-calopt}(Saturation for FSRO with optimal qualification.)
Let $\{g_\alpha\}_{\alpha\in (0,\alpha_0)}$ be a SRM satisfying
hypothesis \textit{H4} and having optimal qualification $\rho$,
$r_\alpha(\lambda)= 1-\lambda g_\alpha(\lambda)$,
$R_\alpha=g_\alpha(T^\ast T)T^\ast$ and suppose that $s_\rho\in
\mathcal{S}$ (where  $s_\rho$ is as defined in (\ref{s-rho})\,).
Assume further that the following hypotheses hold:

\textbf{a)} The function $\rho$ is of local upper type $\beta$,
for some $\beta \geq 0$.

\textbf{b)} There exist positive constants
$\gamma_1,\gamma_2,\lambda^\ast, c_1$, with $\lambda^\ast \leq
\norm{T}^2$ and $c_1>1$ such that

\quad \textbf{i)} $0 \leq r_\alpha(\lambda)\leq 1$, for $\alpha>0$,
$0\leq \lambda\leq \lambda^\ast$;

\quad \textbf{ii)} $r_\alpha(\lambda)\geq\gamma_1$, for
$0\leq\lambda<h(\alpha)\leq \lambda^\ast$, $\alpha\in (0,\alpha_0)$
where $h$ is as in (\ref{eq:4.50}). (Note that by virtue of Theorem
\ref{teo:unicaff} and the fact that $s_\rho\in \mathcal{S}$, there
exists only one function $s \in \mathcal{S}$ satisfying
(\ref{eq:4.50}), that is, $s=s_\rho$.)

\quad \textbf{iii)} $\abs{r_\alpha(\lambda)}$ is non-decreasing with
respect to $\alpha$ for each $\lambda\in (0,\norm{T}^2]$;

\quad \textbf{iv)} $g_\alpha(c_1\alpha)\geq \frac{\gamma_2}{\alpha}$
for $0<c_1\alpha\leq \lambda^\ast$ and

\quad \textbf{v)} $g_\alpha(\lambda)\geq g_\alpha(\tilde{\lambda})$,
for $0<\alpha\leq\lambda\leq \tilde{\lambda}\leq \lambda^\ast$.

\textbf{c)} There exist $\{\lambda_n\}_{n=1}^\infty\subset
\sigma(TT^\ast)$ and $c\geq 1$ such that $\lambda_n \downarrow 0$
and $\frac{\lambda_n}{\lambda_{n+1}}\leq c$ for every $n\in \N$.

Let $\Theta(t)\doteq \sqrt{t}\,\rho(t)$ for $t>0$ and
$X^{s_\rho}\doteq \mathcal{R}(s_\rho(T^\ast T))\setminus \{0\}$.
Then $\psi(x,\delta)\doteq \rho \circ \Theta^{-1}(\delta)$ for $x
\in X^{s_\rho}$ and $\delta \in (0,\Theta (\alpha_0))$, is
saturation function of $\{R_\alpha\}_{\alpha \in (0,\alpha_0)}$ on
$X^{s_\rho}$.
\end{thm}

\smallskip

To prove this theorem we will need three previous lemmas. The
first one is a somewhat technical result, the second one deals
with the existence of an \textit{a-priori} parameter choice rule
leading to a worst total error having an appropriate order of
convergence, while the third one is a converse result.

\medskip
\begin{lem}\label{lem:r-invertible}
Let $\{g_\alpha\}_{\alpha\in (0,\alpha_0)}$ be a SRM, $(\rho,s)$
an order-source pair for $\{g_\alpha\}$ and suppose hypothesis
\textbf{\textit{b.ii)}} of Theorem \ref{teo:sat-calopt} holds. Then
for every $\alpha \in (0,\alpha_0)$ the operator $r_\alpha(T^\ast
T)$ is invertible.
\end{lem}
\smallskip

\begin{proof} Let $\{E_\lambda\}$ be the spectral family of $T^\ast
T$. It suffices to show that for every $\alpha \in (0,\alpha_0)$, $x
\in X$, the function $r_\alpha^{-2}(\lambda)$ is integrable with
respect to the measure $d\norm{E_\lambda x}^2$. Let $\alpha \in
(0,\alpha_0)$ fixed. Since $(\rho,s)$ is an order-source pair for
$\{g_\alpha\}$, there exist  a constant $\gamma>0$ and a function
$h:(0,\alpha_0)\rightarrow \RE^+$ with $\underset{\alpha
\rightarrow0^+}{\lim}h(\alpha)=0$ such that
\begin{equation*}
\frac{s(\lambda)\abs{r_\alpha(\lambda)}}{\rho(\alpha)}\geq \gamma
\quad \forall\; \lambda \in [h(\alpha),+\infty).
\end{equation*}
Therefore
\begin{equation}\label{eq:int1}
\int_{h(\alpha)}^{\norm{T}^2+} \frac 1
{r_\alpha^2(\lambda)}\;d\norm{E_\lambda x}^2\leq
\frac{1}{\gamma^2\rho^2(\alpha)} \int_{h(\alpha)}^{\norm{T}^2+}
s^2(\lambda)\;d\norm{E_\lambda x}^2\leq \frac{\norm{s(T^\ast
T)x}^2}{\gamma^2\rho^2(\alpha)}<+\infty.
\end{equation}

Now, since $\alpha \in (0,\alpha_0)$, it follows from hypothesis
\textbf{\textit{b.ii)}} of Theorem \ref{teo:sat-calopt} that
$r_\alpha(\lambda)\geq \gamma_1>0$ for every  $\lambda \in [0,
h(\alpha))$. Then
\begin{equation}\label{eq:int2}
\int_0^{h(\alpha)} \frac 1 {r_\alpha^2(\lambda)}\;d\norm{E_\lambda
    x}^2\leq  \frac{\norm{x}^2}{\gamma_1^2}<+\infty.
\end{equation}

From (\ref{eq:int1}) and (\ref{eq:int2}) it follows that
$\int_0^{\norm{T}^2+} r_\alpha^{-2}(\lambda)\;d\norm{E_\lambda
    x}^2<+\infty.$ Hence $r_\alpha(T^\ast T)$ is invertible.\hfill
\end{proof}
\medskip
\begin{lem}\label{lem:o-O}
Let $\{g_\alpha\}$ be a SRM, $R_\alpha=g_\alpha(T^\ast T)T^\ast$,
$(\rho,s)$ an order-source pair for $\{g_\alpha\}$ and assume that
hypotheses \textbf{\textit{b.ii)}} and \textbf{\textit{b.iii)}} of
Theorem \ref{teo:sat-calopt} hold. Let
$\varphi:(0,+\infty)\rightarrow\RE^+$ be a continuous,
non-decreasing function satisfying $\underset{\delta \to
0^+}{\lim}\varphi(\delta)=0$ and $x^\ast \in X$, $x^\ast \neq 0$.

\begin{description}
  \item[\textit{I)}] If $\et(x^\ast,\delta)=o(\varphi(\delta))$ for
$\delta \rightarrow 0^+$, then there exists an \textit{a-priori}
parameter choice rule $\tilde{\alpha}(\delta)$ such that
\begin{equation*}
\underset{y^\delta \in \overline{B_\delta(Tx^\ast)}}{\sup}
\norm{R_{\tilde{\alpha}(\delta)}
y^\delta-x^\ast}=o(\varphi(\delta))\quad \textrm{ for } \delta
\rightarrow 0^+.
\end{equation*}
  \item[\textit{II)}] Part \textbf{I)} remains true with $o(\varphi(\delta))$ replaced by $O(\varphi(\delta))$, that is, if $\et(x^\ast,\delta)=O(\varphi(\delta))$ for
$\delta \rightarrow 0^+$, then there exists an \textit{a-priori}
parameter choice rule $\tilde{\alpha}(\delta)$ such that
\begin{equation*}
\underset{y^\delta \in \overline{B_\delta(Tx^\ast)}}{\sup}
\norm{R_{\tilde{\alpha}(\delta)}
y^\delta-x^\ast}=O(\varphi(\delta))\quad \textrm{ for } \delta
\rightarrow 0^+.
\end{equation*}
\end{description}
\end{lem}

\smallskip
\begin{proof}
Let $\varphi$ and $x^\ast \in X$ be as in the hypotheses and suppose
that $\et(x^\ast,\delta)=o(\varphi(\delta))$ for $\delta \rightarrow
0^+$. Then by definition of $\et,$
\begin{equation} \label{eq:11}
\underset{\delta \rightarrow 0^+}{\lim}\frac{\underset{\alpha \in
(0,\alpha_0)}{\inf}\;\underset{y^\delta \in
\overline{B_\delta(Tx^\ast)}}{\sup} \norm{R_\alpha
y^\delta-x^\ast}}{\varphi(\delta)}=\underset{\delta \rightarrow
0^+}{\lim}\,\underset{\alpha \in
(0,\alpha_0)}{\inf}\frac{\underset{y^\delta \in
\overline{B_\delta(Tx^\ast)}}{\sup} \norm{R_\alpha
y^\delta-x^\ast}}{\varphi(\delta)}=0.
\end{equation}
For the sake of simplicity we define:
\begin{equation*}
f(\alpha,\delta)\doteq\frac{\underset{y^\delta \in
\overline{B_\delta(Tx^\ast)}}{\sup} \norm{R_\alpha
y^\delta-x^\ast}}{\varphi(\delta)}\quad \textrm{and}\quad
q(\delta)\doteq \underset{\alpha \in (0,\alpha_0)}{\inf}
f(\alpha,\delta),
\end{equation*}
so that with this notation (\ref{eq:11}) can be written simply as
$\underset{\delta \rightarrow 0^+}{\lim}q(\delta)=0$ and the
objective is to prove the existence of an a-priori parameter
choice rule $\tilde{\alpha}(\delta)$ such that
$\displaystyle\lim_{\delta \to
0^+}f(\tilde{\alpha}(\delta),\delta)=0$. It can be easily proved
that if for certain $\delta_0>0$, $\et(x^\ast, \delta_0)=0$ then
$T^\dag=0.$ Hence $q(\delta)>0$ for every $\delta \in(0,+\infty).$
Also, $q(\delta)$ is continuous for $\delta \in(0,+\infty)$ since
both $\et(x^\ast, \delta)$ and $\varphi(\delta)$ are continuous.
Next, for $n \in \N$ we define
\begin{equation*}
\delta_n\doteq \sup \left\{d>0: q(\delta)\leq \frac 1 n \;\forall\;
\delta\in (0,d) \right\}.
\end{equation*}
Clearly, $\delta_n \downarrow 0$ and since $q$ is continuous for
every $n \in \N$ and every $\delta \in (0, \delta_n]$,
$q(\delta)=\underset{\alpha \in (0,\alpha_0)}{\inf}
f(\alpha,\delta)\leq \frac 1 n$. Then, there exists
$\alpha_n=\alpha_n(\delta_n)\in (0,\alpha_0)$ such that for all $n
\in \N$
\begin{equation}\label{eq:desig}
f(\alpha_n,\delta)\leq \frac 2 n \quad\forall\; \delta \in (0,
\delta_n].
\end{equation}
Now, since $\{\alpha_n\}\subset (0,\alpha_0)$ is a bounded
sequence, there exist $\alpha^\ast \in [0,\alpha_0]$ and
$\{\alpha_{n_k}\}\subset \{\alpha_n\}$ such that $\displaystyle
\lim_{k\to +\infty}\alpha_{n_k}=\alpha^\ast.$ We now define
$\tilde{\alpha}(\delta)\doteq \alpha_{n_k}$ for $\delta \in
(\delta_{n_{k+1}},\delta_{n_k}],\, k=1, 2, ...,$ and
$\tilde{\alpha}(\delta)=\tilde{\alpha}(\delta_{n_1})$ for
$\delta>\delta_{n_1}.$ Then
\begin{equation}\label{eq:lim-alpha}
\lim_{\delta\to 0^+}\tilde{\alpha}(\delta)=\alpha^\ast
\end{equation}
and
\begin{equation*}
0\leq \limsup_{\delta \to 0^+}f(\tilde{\alpha}(\delta),\delta)\leq
\limsup_{k\to
+\infty}\left[\sup_{\delta\in(0,\delta_{n_k}]}f(\alpha_{n_k},\delta)
\right]\leq \limsup_{k\to +\infty}\frac 2 {n_k}=0,
\end{equation*}
where the last inequality follows from (\ref{eq:desig}).

Hence,
\begin{equation}\label{eq:lim-f}
\lim_{\delta\to 0^+}f(\tilde{\alpha}(\delta), \delta)=0.
\end{equation}

It remains to be shown that $\tilde{\alpha}(\delta)$ is an
admissible parameter choice rule, for which it suffices to prove
that $\displaystyle\lim_{\delta\to 0^+}\tilde{\alpha}(\delta)=0$,
i.e. that $\alpha^\ast=0$. If $\alpha^\ast>0$, it follows from
(\ref{eq:lim-alpha}) that there exists $\delta_0>0$ such that
$\tilde{\alpha}(\delta)>\frac{\alpha^\ast}{2}$ for all $\delta \in
(0,\delta_0)$. Hypothesis \textbf{\textit{b.iii)}} of Theorem
\ref{teo:sat-calopt} then implies that for every $\delta \in
(0,\delta_0)$, $\abs{r_{\tilde{\alpha}(\delta)}(\lambda)}\geq
\abs{r_{\frac{\alpha^\ast}{2}}(\lambda)}$ for all $\lambda \in (0,
\norm{T}^2].$ Therefore for every $\delta \in (0,\delta_0)$,
\begin{eqnarray}\nonumber
  \norm{r_{\tilde{\alpha}(\delta)}(T^\ast
T)x^\ast}^2 &=&
\int_0^{\norm{T}^2+}r_{\tilde{\alpha}(\delta)}^2(\lambda)\;d\norm{E_\lambda
    x^\ast}^2 \\ \nonumber
&\geq&\int_0^{\norm{T}^2+}r_{\frac{\alpha^\ast}{2}}^2(\lambda)\;d\norm{E_\lambda
    x^\ast}^2 \\ \label{eq:des1}
    &=&\norm{r_{\frac{\alpha^\ast}{2}}(T^\ast
T)x^\ast}^2.
\end{eqnarray}
Now, for all $\delta \in (0,\delta_0)$,
\begin{eqnarray*}
\underset{y^\delta \in \overline{B_\delta(Tx^\ast)}}{\sup}
\norm{R_{\tilde{\alpha}(\delta)} y^\delta-x^\ast} &\geq&
\norm{R_{\tilde{\alpha}(\delta)} Tx^\ast-x^\ast}
   = \norm{\left(I-g_{\tilde{\alpha}(\delta)}(T^\ast T)T^\ast T\right)x^\ast} \\
   &=& \norm{r_{\tilde{\alpha}(\delta)}(T^\ast T)x^\ast}
   \geq\norm{r_{\frac{\alpha^\ast}{2}}(T^\ast T)x^\ast},
\end{eqnarray*}
where the last inequality follows from (\ref{eq:des1}). Dividing
through by $\varphi(\delta)$, taking limit for $\delta \rightarrow
0^+$, and using the definition of $f(\alpha,\delta)$ and
(\ref{eq:lim-f}) we conclude that
$\norm{r_{\frac{\alpha^\ast}{2}}(T^\ast T)x^\ast}=0.$ Now since
$\frac{\alpha^\ast}{2}<\alpha_0$, $(\rho,s)$ is an order-source pair
for $\{g_\alpha\}$ and hypothesis \textbf{\textit{b.ii)}} of Theorem
\ref{teo:sat-calopt} holds, it follows from Lemma
\ref{lem:r-invertible} that $r_{\frac{\alpha^\ast}{2}}(T^\ast T)$ is
invertible. Therefore $x^\ast=0$, contradicting the hypothesis that
$x^\ast\neq 0$. Hence, $\alpha^\ast$ must be equal to zero, as
wanted.

We proceed now to prove the second part of the Lemma. Suppose that
there exists $x^\ast \in X$, $x^\ast\neq 0$ such that
$\et(x^\ast,\delta)=O(\varphi(\delta))$ as $\delta \rightarrow 0^+$.
Then there exist positive constants $k$ and $d$ such that
$\underset{\alpha \in (0,\alpha_0)}{\inf}f(\alpha, \delta)\leq k$
for every $\delta \in (0,d)$, where $f(\alpha,\delta)$ is as
previously defined. Let $\{\delta_n\}_{n\in \N}\subset (0,d)$ be
such that $\delta_n\downarrow 0$ and $\alpha_n=\alpha_n(\delta_n)\in
(0,\alpha_0)$ such that
\begin{equation*}
    f(\alpha_n,\delta)\leq k+\delta_n, \;\forall\, \delta \in
    (0,d),\;\forall \,n \in \N.
\end{equation*}
Without loss of generality we assume that the sequence
$\{\alpha_n\}$ converges (since if that is not the case, we can take
a subsequence which does). Now, like in the previously case, by
defining $\tilde{\alpha}(\delta)=\alpha_n$ for $\delta \in
(\delta_{n+1},\delta_n], n=1,2,...,$ and
$\tilde{\alpha}(\delta)=\alpha(\delta_1)$ for $\delta >\delta_1$,
since $\delta_n \downarrow 0$ it follows that
$f(\tilde{\alpha}(\delta),\delta)\leq k+\delta_1$ for every $\delta
\in (0,d)$ and therefore
\begin{equation*}
\underset{y^\delta \in \overline{B_\delta(Tx^\ast)}}{\sup}
\norm{R_{\tilde{\alpha}(\delta)}
y^\delta-x^\ast}=O(\varphi(\delta))\quad \textrm{as } \delta
\rightarrow 0^+.
\end{equation*}

Following the same steps as in the proof of Part I we obtain that
$\displaystyle \lim_{\delta \to 0^+}\tilde{\alpha}(\delta)=0$, i.e.
$\tilde{\alpha}(\delta)$ is an admissible parameter choice rule.
This concludes the proof of the lemma.\hfill
\end{proof}
%
\smallskip
\begin{lem}\label{lem:sup-inf} Let $\{g_\alpha\}_{\alpha\in (0,\alpha_0)}$ be a
SRM, $r_\alpha(\lambda)\doteq 1-\lambda g_\alpha(\lambda)$,
$R_\alpha=g_\alpha(T^\ast T)T^\ast$, $(\rho, s)$ an order-source
pair for $\{g_\alpha\}$, $\Theta(t)\doteq \sqrt{t}\,\rho(t)$ for
$t>0$, and suppose that:

\textbf{a)} The function $\rho$ is of local upper type $\beta$,
for some $\beta \geq 0$.

\textbf{b)} There exist positive constants
$\gamma_1,\gamma_2,\lambda^\ast, c_1$, with $\lambda^\ast \leq
\norm{T}^2$ and $c_1>1$ such that

\quad \textbf{i)} $0 \leq r_\alpha(\lambda)\leq 1$, for
$\alpha>0$, $0\leq \lambda\leq \lambda^\ast$;

\quad \textbf{ii)} $r_\alpha(\lambda)\geq\gamma_1$, for
$0\leq\lambda<h(\alpha)\leq \lambda^\ast$, $\alpha\in (0,\alpha_0)$
where $h$ is as in (\ref{eq:4.50});

\quad \textbf{iii)} $\abs{r_\alpha(\lambda)}$ is non-decreasing with
respect to $\alpha$ for each $\lambda\in (0,\norm{T}^2]$;

\quad \textbf{iv)} $g_\alpha(c_1\alpha)\geq \frac{\gamma_2}{\alpha}$
for $0<c_1\alpha\leq \lambda^\ast$ and

\quad \textbf{v)} $g_\alpha(\lambda)\geq g_\alpha(\tilde{\lambda})$,
for $0<\alpha\leq\lambda\leq \tilde{\lambda}\leq \lambda^\ast$.

\textbf{c)} There exist $\{\lambda_n\}_{n=1}^\infty\subset
\sigma(TT^\ast)$ and $c\geq 1$ such that $\lambda_n \downarrow 0$
and $\frac{\lambda_n}{\lambda_{n+1}}\leq c$ for every $n\in \N$.

If for some $x \in X$ we have that
\begin{equation}\label{eq:19} \underset{y^\delta \in
\overline{B_\delta(Tx)}}{\sup}\;\underset{\alpha \in
(0,\alpha_0)}{\inf} \norm{R_{\alpha}
y^\delta-x}=O(\rho(\Theta^{-1}(\delta))) \quad \textrm{when } \delta
\rightarrow 0^+,
\end{equation}
then $x \in \mathcal{R}(s(T^\ast T))$. In particular, if $\rho$ is
optimal qualification of $\{g_\alpha\}$ and $s_\rho\in
\mathcal{S}$, then $x \in \mathcal{R}(s_\rho(T^\ast T))$.
\end{lem}
\medskip

\begin{proof} Without loss of generality we may assume that
$\alpha_0\le\frac{\lambda^\ast}{c_1}$ and $x \neq 0$ (this is so
because hypotheses \textbf{\textit{a)}} and \textit{\textbf{c)}}
are independent of $\alpha_0$ and if  \textit{\textbf{b)}} holds
for $\alpha \in (0,\alpha_0)$ then it holds for $\alpha\in
(0,\hat{\alpha}_0)$ for every $\hat{\alpha}_0<\alpha_0$ with the
same constants, while if $x=0$ the result of the Lemma is
trivial).

Let $\bar{\alpha}\in \sigma(TT^\ast)$ be such that
$0<c_1\,\bar{\alpha}\leq \alpha_0$ (hypothesis
\textit{\textbf{c)}} guarantees the existence of such
$\bar{\alpha}$), and define
\begin{equation*}
    \bar{\delta}=\bar{\delta}(\bar\alpha)\doteq
    \frac{\bar{\alpha}^{1/2}}{\gamma_2}\norm{R_{\bar{\alpha}}Tx-x}=\frac{\bar{\alpha}^{1/2}}{\gamma_2}\norm{r_{\bar{\alpha}}(T^\ast T)x}.
\end{equation*}
Then, clearly the equation
\begin{equation}\label{eq:14}
\norm{R_\alpha Tx -x}^2=\frac{(\gamma_2\,\bar{\delta})^2}{\alpha}
\end{equation}
in the unknown $\alpha$, has $\alpha=\bar{\alpha}$ as a solution.
Moreover, since $\norm{R_\alpha Tx-x}^2=\int_0^{\|T\|^{2\,+}}
r_\alpha^2(\lambda)\, d\norm{E_\lambda x}^2$ and $x \neq 0$,
hypotheses \textbf{\textit{b.ii)}} and \textbf{\textit{b.iii)}}
imply that the function $\mu(\alpha)\doteq\alpha\norm{R_\alpha
Tx-x}^2$ is strictly increasing for $\alpha$ in $(0,\alpha_0)$.
Hence, $\alpha =\bar{\alpha}\doteq\eta(\bar\delta)$ (where
$\eta(\delta)=\mu^{-1}(\delta)$\,) is the unique solution of
(\ref{eq:14}). Note that if $\bar{\alpha}\rightarrow 0^+$ then
$\bar{\delta}\rightarrow 0^+$. Moreover, by hypothesis
\textbf{\textit{b.iii)}} and Lemma \ref{lem:r-invertible}, it
follows immediately by Fatou's Lemma that $\bar{\delta}\rightarrow
0^+$ only if $\bar{\alpha}\rightarrow 0^+$.

Now, for $\delta>0$ define
\begin{equation}\label{eq:16}
    \bar{y}^{\,\delta}\doteq Tx-\delta G_{\bar{\alpha}}z,
\end{equation}
where $G_{\bar{\alpha}}\doteq
F_{c_1\,\bar{\alpha}}-F_{\bar{\alpha}}$ with $\{F_\lambda\}$ being
the spectral family associated to $TT^\ast$ and
\begin{equation*}
    z\doteq \left\{
\begin{array}{ll}
    \norm{G_{\bar{\alpha}}Tx}^{-1}Tx, & \hbox{if $G_{\bar{\alpha}}Tx\neq
0$,} \\
    \textrm{arbitrary with} \norm{G_{\bar{\alpha}}z}=1, &
    \hbox{in other case.} \\
\end{array}
\right.
\end{equation*}
Note that since $\bar{\alpha}\in \sigma(TT^\ast)$ and $c_1>1$ it
follows that $G_{\bar{\alpha}}$ is not the null operator and
therefore the definition makes sense. Note also that
$\norm{\bar{y}^{\,\delta}-Tx}=\delta$, which implies that
$\bar{y}^{\,\delta} \in \overline{B_\delta(Tx)}$.

Now, by using (\ref{eq:16}) and the fact that $g_\alpha(T^\ast
T)T^\ast= T^\ast g_\alpha(TT^\ast)$ it follows that for every
$\alpha \in (0,\alpha_0)$ and $\delta>0$,
\begin{eqnarray}\label{eq:12} \nonumber
\seq{R_\alpha Tx-x\right.\left.,R_\alpha (\bar{y}^{\,\delta}-Tx)}
&=& \seq{g_\alpha(T^\ast T)T^\ast Tx-x
,-g_\alpha(T^\ast T)T^\ast\,\delta G_{\bar{\alpha}}z} \\
\nonumber  &=& \delta \seq{g_\alpha(T^\ast T)T^\ast Tx-x,-T^\ast
g_\alpha(TT^\ast)G_{\bar{\alpha}}z} \\ \nonumber &=& \delta
\seq{Tg_\alpha(T^\ast T)T^\ast
Tx-Tx,-g_\alpha(TT^\ast)G_{\bar{\alpha}}z} \\ \nonumber  &=&
\delta \seq{(TT^\ast
g_\alpha(TT^\ast)-I)Tx,-g_\alpha(TT^\ast)G_{\bar{\alpha}}z}
\\ \nonumber
&=& \delta
\seq{-r_\alpha(TT^\ast)Tx,-g_\alpha(TT^\ast)G_{\bar{\alpha}}z}\\
&=& \delta \int_0^{\norm{T}^2+}
   r_\alpha(\lambda)g_\alpha(\lambda)\,d\seq{F_\lambda
   Tx,G_{\bar{\alpha}}z}.
\end{eqnarray}

Now by hypothesis \textbf{\textit{b.i)}} and since
$c_1\bar{\alpha}\leq \lambda^\ast$ one has that both
$g_\alpha(\lambda)$ and $r_\alpha(\lambda)$ are nonnegative for
all $\lambda \in [0,c_1\bar{\alpha}]$. Also, from the definitions
of $G_{\bar{\alpha}}$ and $z$ it follows immediately that the
function $m(\lambda)\doteq \seq{F_\lambda Tx,G_{\bar{\alpha}}z}$
for $\lambda\in[0,c_1\bar\alpha]$ is real and non-decreasing and
therefore
\begin{equation}\label{eq:28}
    \int_0^{c_1\bar{\alpha}}r_\alpha(\lambda)g_\alpha(\lambda)\,d\seq{F_\lambda
   Tx,G_{\bar{\alpha}}z}\geq 0.
\end{equation}
On the other hand, since $m(\lambda)= \seq{Tx,F_\lambda
G_{\bar{\alpha}}z}$ and $F_\lambda
G_{\bar{\alpha}}=G_{\bar{\alpha}}$ for every $\lambda \geq
c_1\bar{\alpha}$, it follows that $m(\lambda)$ is constant for $\lambda \geq c_1\bar{\alpha}$ and therefore
\begin{equation}\label{eq:29}
 \int_{c_1\bar{\alpha}}^{\norm{T}^2+}r_\alpha(\lambda)g_\alpha(\lambda)\,d\seq{F_\lambda
   Tx,G_{\bar{\alpha}}z}= 0.
\end{equation}
From (\ref{eq:28}) and (\ref{eq:29}) we conclude that
\begin{equation*}
\int_0^{\norm{T}^2+}
   r_\alpha(\lambda)g_\alpha(\lambda)\,d\seq{F_\lambda
   Tx,G_{\bar{\alpha}}z}\geq0,
\end{equation*}
which, by virtue of (\ref{eq:12}), implies that
\begin{equation}\label{eq:30}
\seq{R_\alpha Tx-x,R_\alpha (\bar{y}^{\,\delta}-Tx)}\geq 0.
\end{equation}

Hence, for every $\alpha\in (0,\alpha_0)$, $\delta>0$ and
$\bar{\alpha}\in \sigma(TT^\ast)$ such that $c_1\bar{\alpha}\leq
\lambda^\ast$ we obtain the following estimate:
\begin{eqnarray}\label{eq:20}\nonumber
\norm{R_\alpha \bar{y}^{\,\delta}-x}^2 &=& \norm{R_\alpha Tx-x}^2
+\norm{R_\alpha (\bar{y}^{\,\delta}-Tx)}^2 + 2\seq{R_\alpha
Tx-x,R_\alpha (\bar{y}^{\,\delta}-Tx)}\\ \nonumber    &=&
\norm{R_\alpha Tx-x}^2 + \delta^2 \norm{g_\alpha(T^\ast T)T^\ast
G_{\bar{\alpha}}z}^2 + 2\seq{R_\alpha Tx-x,R_\alpha
(\bar{y}^{\,\delta}-Tx)}\text{(by (\ref{eq:16})\,)}\\
\nonumber &\geq&\norm{R_\alpha Tx-x}^2
+ \delta^2 \norm{T^\ast g_\alpha( TT^\ast)G_{\bar{\alpha}}z}^2 \quad\qquad\qquad\qquad\qquad\qquad\qquad\;\;\text{(by (\ref{eq:30})\,)}  \\
\nonumber &=&\norm{R_\alpha Tx-x}^2
+ \delta^2 \norm{(TT^\ast)^{\frac12} g_\alpha(TT^\ast)G_{\bar{\alpha}}z}^2 \\
\nonumber  &=& \norm{R_\alpha Tx-x}^2 + \delta^2
\int_0^{{\|T\|^2}\,^+}\lambda\,
g_{\alpha}^2(\lambda)\, d\norm{F_\lambda G_{\bar{\alpha}}z}^2 \\
&\geq& \norm{R_\alpha Tx-x}^2 + \delta^2
\int_{\bar{\alpha}}^{c_1\,\bar{\alpha}}\lambda
\,g_{\alpha}^2(\lambda)\, d\norm{F_\lambda G_{\bar{\alpha}}z}^2.
\end{eqnarray}

We now consider two different cases for $\alpha\in (0,\alpha_0)$.

\underline{Case I}: $\alpha \leq \bar{\alpha}$. Since
$c_1\bar{\alpha}\leq\lambda^\ast$ and $c_1>1$, it follows from
hypothesis \textbf{\textit{b.v)}} that
\begin{equation}\label{eq:17}
g_\alpha(\lambda)\geq g_\alpha(c_1\bar{\alpha})\geq
g_\alpha(\lambda^\ast) \quad\textrm{for every}\quad \lambda \in
[\bar{\alpha},c_1\,\bar{\alpha}].
\end{equation}
On the other hand, from hypothesis \textbf{\textit{b.i)}} it follows
that $r_\alpha(\lambda^\ast)\leq 1$, which implies that
$\lambda^\ast\,g_\alpha(\lambda^\ast)\geq 0$ and therefore,
$g_\alpha(\lambda^\ast)\geq 0.$ It then follows from (\ref{eq:17})
that $g_\alpha^2(\lambda)\geq g_\alpha^2(c_1\,\bar{\alpha})$ for
every $\lambda \in [\bar{\alpha},c_1\,\bar{\alpha}].$ Then,
\begin{eqnarray} \label{eq:21} \nonumber
    \int_{\bar{\alpha}}^{c_1\,\bar{\alpha}}\lambda
\,g_{\alpha}^2(\lambda)\, d\norm{F_\lambda
    G_{\bar{\alpha}}z}^2 &\geq& \bar{\alpha}
\,g_{\alpha}^2(c_1\,\bar{\alpha})\int_{\bar{\alpha}}^{c_1\,\bar{\alpha}}
d\norm{F_\lambda
    G_{\bar{\alpha}}z}^2 \\ \nonumber
    &=& \bar{\alpha}
\,g_{\alpha}^2(c_1\,\bar{\alpha})\left(\norm{F_{c_1\,\bar{\alpha}}
    G_{\bar{\alpha}}z}^2-\norm{F_{\bar{\alpha}}
    G_{\bar{\alpha}}z}^2\right) \\ \nonumber
    &=& \bar{\alpha}
\,g_{\alpha}^2(c_1\,\bar{\alpha})\norm{G_{\bar{\alpha}}z}^2 \\
    &=& \bar{\alpha}\,g_{\alpha}^2(c_1\,\bar{\alpha}),
\end{eqnarray}
where the second to last equality follows from the definition of
    $G_{\bar{\alpha}}$ and the spectral property $F_\lambda
    F_\mu=F_{\min\{\lambda,\mu\}}.$

At the same time, the hypotheses \textbf{\textit{b.i)}} and
\textbf{\textit{b.iii)}} imply that $g_\alpha(\lambda)$ is
non-increasing as a function of $\alpha$ for each fixed $\lambda \in
[0,\lambda^\ast]$. Since $\alpha\leq \bar{\alpha}$ and
$c_1\,\bar{\alpha}\leq \lambda^\ast$, we then have that
\begin{equation}\label{eq:18} g_\alpha(c_1\,\bar{\alpha})\geq
g_{\bar{\alpha}}(c_1\,\bar{\alpha}),
\end{equation}
and from hypothesis \textbf{\textit{b.iv)}} we also have that
\begin{equation}\label{eq:18-1}
g_{\bar{\alpha}}(c_1\,\bar{\alpha})\geq \frac{\gamma_2}{\bar{\alpha}}>0.
\end{equation}
From (\ref{eq:18}) and (\ref{eq:18-1}) we conclude that
\begin{equation}\label{eq:18-2}
g_\alpha^2(c_1\,\bar{\alpha})\geq
\left(\frac{\gamma_2}{\bar{\alpha}}\right)^2.
\end{equation}

Substituting (\ref{eq:18-2}) into (\ref{eq:21}) we obtain
$\int_{\bar{\alpha}}^{c_1\,\bar{\alpha}}\lambda
\,g_{\alpha}^2(\lambda)\, d\norm{F_\lambda
    G_{\bar{\alpha}}z}^2 \geq \frac{\gamma_2^2}{\bar{\alpha}}$,
which, by virtue of (\ref{eq:20}) implies that if $\alpha\leq
\bar{\alpha}$, then $\norm{R_\alpha \bar{y}^{\,\delta}-x}^2\geq
\frac{(\gamma_2\,\delta)^2}{\bar{\alpha}}$.

\smallskip
\underline{Case II}: $\alpha>\bar{\alpha}$. In this case, it
follows from hypothesis \textbf{\textit{b.iii)}} that
$r_\alpha^2(\lambda)\ge r_{\bar{\alpha}}^2(\lambda)$ for every
$\lambda \in (0,\norm{T}^2]$. Then,
\begin{equation*}
\norm{R_\alpha Tx-x}^2=\int_0^{\|T\|^{2\,+}} r_\alpha^2(\lambda)\,
d\norm{E_\lambda x}^2 \ge \int_0^{\|T\|^{2\,+}}
r_{\bar{\alpha}}^2(\lambda)\, d\norm{E_\lambda
x}^2=\norm{R_{\bar{\alpha}} Tx-x}^2,
\end{equation*}
which, together with (\ref{eq:20}) imply that $\norm{R_\alpha
\bar{y}^{\,\delta}-x}^2\ge \norm{R_{\bar{\alpha}} Tx-x}^2$.

\smallskip
Summarizing the results of cases I and II, we obtain that for
every $\alpha\in(0,\alpha_0)$, $\delta>0$,
$\bar\alpha\in\sigma(TT^\ast)$ with $c_1\bar\alpha\le\alpha_0$ and
$\bar y^\delta$ as in (\ref{eq:16}), there holds:
\begin{eqnarray} \label{eq:18-3} \nonumber
\norm{R_\alpha \bar{y}^{\,\delta}-x}^2&\geq& \left\{
\begin{array}{ll}
\frac{(\gamma_2\,\delta)^2}{\bar{\alpha}}, & \hbox{if $0<\alpha
\leq
\bar{\alpha}$,} \\
    \norm{R_{\bar{\alpha}} Tx-x}^2, & \hbox{if $\bar\alpha<\alpha<\alpha_0$} \\
\end{array}
\right. \\
&\geq&\min\left\{\norm{R_{\bar{\alpha}} Tx-x}^2,
\frac{(\gamma_2\,\delta)^2}{\bar{\alpha}}\right\}.
\end{eqnarray}
Then
\begin{eqnarray*}
\min\left\{\norm{R_{\bar{\alpha}} Tx-x}, \frac{\gamma_2
\,\delta}{\sqrt{\bar{\alpha}}}\right\}&=&
\left(\min\left\{\norm{R_{\bar{\alpha}} Tx-x}^2,
\frac{(\gamma_2\,\delta)^2}{\bar{\alpha}}\right\}\right)^{1/2} \\
&\leq& \underset{\alpha \in (0,\alpha_0)}{\inf} \norm{R_\alpha
\bar{y}^{\,\delta}-x}\quad \quad \qquad \qquad \qquad
\parbox{2.5cm}{(by (\ref{eq:18-3}))}  \\
&\leq& \underset{y^\delta \in \overline{B_\delta(Tx)}}{\sup}\; \underset{\alpha \in (0,\alpha_0)}{\inf} \norm{R_\alpha y^\delta-x}\quad \; \parbox{4cm}{(since $\bar{y}^{\,\delta} \in \overline{B_\delta(Tx)}$)} \\
&=& O(\rho(\Theta^{-1}(\delta))) \quad \textrm{for} \quad \delta
\rightarrow
 0^+ \quad \parbox{4cm}{(by hypothesis),}
\end{eqnarray*}
and since $\bar{\alpha}=\eta(\bar{\delta})$ solves equation
(\ref{eq:14}), the previous inequality implies that
\begin{equation}\label{eq:22}
\norm{R_{\eta(\bar{\delta})} Tx-x}=\frac{\gamma_2
\,\bar{\delta}}{\sqrt{\bar{\alpha}}}=O(\rho(\Theta^{-1}(\bar{\delta})))
\quad \textrm{for} \quad \bar{\delta} \rightarrow 0^+,
\end{equation}
and therefore
\begin{equation}\label{eq:23}
\frac{\bar{\delta}}{\rho(\Theta^{-1}(\bar{\delta}))}=O\left(\sqrt{\eta(\bar{\delta})}\right)
\quad \textrm{for} \quad \bar{\delta} \rightarrow 0^+.
\end{equation}
Now, since for every $\delta
>0$ one has $\delta=\Theta(\Theta^{-1}(\delta))$, it follows from the
definition of $\Theta$ that
$\delta=\sqrt{\Theta^{-1}(\delta)}\,\rho(\Theta^{-1}(\delta))$.
Then, from (\ref{eq:23}) we obtain that
$\sqrt{\Theta^{-1}(\bar{\delta})}=O(\sqrt{\eta(\bar{\delta})})$
for $\bar{\delta} \rightarrow 0^+.$ From this and (\ref{eq:22}) we
then deduce that
\begin{equation}\label{eq:27}
\norm{R_{\eta(\bar{\delta})}Tx-x}=O(\rho(\eta(\bar{\delta})))
\textrm{\;\, for\;\,} \bar{\delta} \rightarrow 0^+, \; \bar\delta=
\mu(\bar{\alpha}),\;\bar\alpha \in \sigma(TT^\ast),\;
c_1\,\bar\alpha\le\alpha_0.
\end{equation}
\medskip
Now let
$L\doteq\max\left\{\lambda_j\;:\;\lambda_j\le\frac{\alpha_0}{c_1}
\right\}$. Then, since $\lambda_n\downarrow 0$, for any
$\alpha\in(0,L]$ there exist a unique $n=n(\alpha)\in\mathbb{N}$
such that $\lambda_{n+1}<\alpha\le\lambda_n$ (note that
$n(\alpha)\to\infty$ if and only if $\alpha\to 0^+$). Then for
$\alpha \in (0,L]$ and $n=n(\alpha)$ so defined we have that
\begin{eqnarray}\label{eq:24}\nonumber
  \norm{R_\alpha Tx-x}^2 &=& \int_0^{\|T\|^{2\,+}} r_\alpha^2(\lambda)\, d\norm{E_\lambda
  x}^2\\ \nonumber
   &\leq& \int_0^{\|T\|^{2\,+}} r_{\lambda_n}^2(\lambda)\,d\norm{E_\lambda
   x}^2  ,\qquad \text{(by hypothesis \textit{\textbf{b.iii)}}\; )}   \\ \nonumber
   &=&\norm{R_{\lambda_n} Tx-x}^2\\
   &=&O(\rho^2(\lambda_n)),\quad \text{(by virtue of (\ref{eq:27}), with
   $\bar\delta=\mu(\lambda_n)$\,).}
\end{eqnarray}
Also, from hypothesis \textbf{\textit{c)}} we have that
$\lambda_n\leq c\,\lambda_{n+1}$ for all $n\in\mathbb{N}$, and
since $\rho$ is non-decreasing and positive (since $\rho \in
\mathcal{O}$) it follows that
\begin{equation}\label{eq:25}
\rho^2(\lambda_n)\leq \rho^2(c\,\lambda_{n+1}),\quad\forall
n\in\mathbb{N}.
\end{equation}
Now since $c\geq 1$ and  by hypothesis \textbf{\textit{a)}} $\rho$
is of local upper type $\beta$, there exists a positive constant
$d$ such that
\begin{equation}\label{eq:26}
\rho(c\,\lambda_{n+1})\leq d\,c^\beta\rho\left(\frac 1 c \,
c\,\lambda_{n+1}\right)= d\,c^\beta\rho(\lambda_{n+1}),\quad
\forall n\in \mathbb{N}.
\end{equation}
Hence
\begin{eqnarray}\label{eq:rho-lambda}
\rho\left(\lambda_{n(\alpha)}\right)&\le &
\rho\left(c\lambda_{n(\alpha)+1}\right)\quad\qquad\text{(by (\ref{eq:25})\;)}\nonumber\\
&\le & dc^\beta \rho\left(\lambda_{n(\alpha)+1}\right) \qquad \text{(by (\ref{eq:26})\;)}\nonumber\\
&\le & dc^\beta \rho\left(\alpha\right) \qquad \text{(since
$\lambda_{n(\alpha)+1}<\alpha$ and $\rho\in\cal{O}$)}.
\end{eqnarray}

From (\ref{eq:24}) and (\ref{eq:rho-lambda}) it follows that
$\norm{R_\alpha Tx-x}=O(\rho(\alpha))$ for $\alpha \rightarrow
0^+$. Since $T^\dag T$ is the projection on
$\mathcal{N}(T)^\perp=X$ (since $T$ is invertible), we have that
$x= T^\dag T x=T^\dag y$. Then $\norm{(R_\alpha -
T^\dag)y}=O(\rho(\alpha))$ for $\alpha \rightarrow 0^+$. Finally,
Theorem \ref{teo:gen-411} implies that $T^\dag y=x\in
\mathcal{R}(s(T^\ast T))$. This concludes the proof of the
Lemma.\hfill
\end{proof}

\begin{rem}
Note that since $\underset{y^\delta \in
\overline{B_\delta(Tx)}}{\sup}\;\underset{\alpha \in
(0,\alpha_0)}{\inf} \norm{R_{\alpha} y^\delta-x} \;\le\;
\et(x,\delta)$, hypothesis (\ref{eq:19}) of the preceding Lemma
holds if \; $\et(x,\delta)= O(\rho(\Theta^{-1}(\delta)))$ for
$\delta\to 0^+$.
\end{rem}

\bigskip
Having stated and proved the previous three lemmas, we are now
ready to prove Theorem \ref{teo:sat-calopt}.

\smallskip

\textit{Proof of Theorem \ref{teo:sat-calopt}.} As in Lemma
\ref{lem:sup-inf}, without loss of generality we may assume that
$\alpha_0\le\frac{\lambda^\ast}{c_1}$. We will show that
$\psi(x,\delta)\doteq \rho \circ \Theta^{-1}(\delta)$ for $x \in
X^{s_\rho}$ and $\delta \in (0, \Theta(\alpha_0))$, is saturation
function of $\{R_\alpha\}_{\alpha\in (0,\alpha_0)}$ on
$X^{s_\rho}$ (see Definition \ref{def:satur}).

First we note that since $\{g_\alpha\}$ satisfies (H4) and $\rho$
is continuous ($\rho$ being of local upper type), by virtue of
Lemma \ref{lem:cotasup} one has that $\psi \in
\mathcal{U}_{X^{s_\rho}}(\et)$. Next we will show that $\psi$
satisfies the \textit{S1} condition for saturation on
$X^{s_\rho}$. Suppose that it is not true, i.e. that there exist
$x^\ast \in X$, $x^\ast\neq 0$ and $x \in X^{s_\rho}$ such that
$\underset{\delta
\rightarrow0^+}{\limsup}\,\frac{\et(x^\ast,\delta)}{\psi(x,\delta)}=0$.
Then $\et(x^\ast,\delta)=o(\psi(x,\delta))$ as $\delta \rightarrow
0^+$ and from Lemma \ref{lem:o-O} I) it follows that there exists
an \textit{a-priori} parameter choice rule
 $\hat\alpha(\delta)$ such that
\begin{equation}\label{eq:o}
\underset{y^\delta \in \overline{B_\delta(Tx^\ast)}}{\sup}
\norm{R_{\hat\alpha(\delta)}
y^\delta-x^\ast}=o(\psi(x,\delta))\quad \textrm{ for } \delta
\rightarrow 0^+.
\end{equation}

On the other hand, from hypothesis \textit{\textbf{c)}} it follows
that there exists $\bar{\alpha} \in\sigma(TT^\ast)$ such that
$0<c_1\bar{\alpha}\leq \alpha_0$ and $h(\bar{\alpha})<\norm{T}^2$.
Now define
\begin{equation*}
    \bar{\delta}=\bar{\delta}(\bar\alpha)\doteq
    \frac{\bar{\alpha}^{1/2}}{\gamma_2}\norm{R_{\bar{\alpha}}Tx^\ast-x^\ast}=\frac{\bar{\alpha}^{1/2}}{\gamma_2}
    \norm{r_{\bar{\alpha}}(T^\ast T)x^\ast}
\end{equation*}
and for $\delta>0$,
\begin{equation}\label{eq:ydj}
    \bar{y}^{\,\delta}\doteq Tx^\ast-\delta \,G_{\bar{\alpha}}z,
\end{equation}
as in (\ref{eq:16}). Following the same steps as in the proof of Lemma \ref{lem:sup-inf} we obtain as in (\ref{eq:18-3}) that for every $\alpha\in(0,\alpha_0)$, $\delta>0$, $\bar\alpha\in\sigma(TT^\ast)$ with $c_1\bar\alpha\le\alpha_0$ and $h(\bar{\alpha})<\norm{T}^2$, and
$\bar y^\delta$ as in (\ref{eq:ydj}), there holds:
\begin{equation}\label{eq:des2}
\norm{R_\alpha \bar{y}^{\,\delta}-x^\ast}^2\geq\min\left\{\norm{R_{\bar{\alpha}}Tx^\ast-x^\ast}^2,
\frac{(\gamma_2\,\delta)^2}{\bar{\alpha}}\right\}.
\end{equation}
Then for $\delta>0$ such that $\hat\alpha(\delta)\in (0,
\alpha_0),$
\begin{eqnarray*}
\min\left\{\norm{R_{\bar{\alpha}}Tx^\ast-x^\ast}, \frac{\gamma_2
\,\delta}{\sqrt{\bar{\alpha}}}\right\}&=&
\left(\min\left\{\norm{R_{\bar{\alpha}}Tx^\ast-x^\ast}^2,
\frac{(\gamma_2\,\delta)^2}{\bar{\alpha}}\right\}\right)^{1/2} \\
&\leq& \underset{\alpha \in (0,\alpha_0)}{\inf} \norm{R_\alpha
\bar{y}^{\,\delta}-x^\ast}\quad \quad \qquad \qquad \qquad \quad
\parbox{3.5cm}{(by (\ref{eq:des2}))}  \\
&\leq&\norm{R_{\hat\alpha(\delta)}
\bar{y}^{\,\delta}-x^\ast}\qquad \qquad \qquad
\parbox{4cm}{(since $\hat\alpha(\delta)\in (0, \alpha_0)$)}\\
&\leq& \underset{y^\delta \in \overline{B_\delta(Tx^\ast)}}{\sup}\;
\norm{R_{\hat\alpha(\delta)} y^\delta-x^\ast}\quad \quad \parbox{4cm}{(since $\bar{y}^{\,\delta} \in \overline{B_\delta(Tx^\ast)}$)} \\
&=& o(\rho(\Theta^{-1}(\delta)))\quad
\textrm{ for } \delta \rightarrow 0^+ \quad \qquad \qquad \parbox{4cm}{(by (\ref{eq:o})),}
\end{eqnarray*}
and since $\bar{\alpha}=\eta(\bar{\delta})$ solves equation
(\ref{eq:14}) with $x=x^\ast$, the previous inequality implies that
\begin{equation*}
\norm{R_{\eta(\bar{\delta})} Tx^\ast-x^\ast}=\frac{\gamma_2
\,\bar{\delta}}{\sqrt{\bar{\alpha}}}=o(\rho(\Theta^{-1}(\bar{\delta})))
\quad \textrm{for} \quad \bar{\delta} \rightarrow 0^+.
\end{equation*}
Following analogous steps as in the proof of Lemma
\ref{lem:sup-inf} we obtain, as in (\ref{eq:27}), that
\begin{equation}\label{eq:nueva}
\norm{R_{\eta(\bar{\delta})}Tx^\ast-x^\ast}=o(\rho(\eta(\bar{\delta})))
\textrm{\; for\;} \bar{\delta} \rightarrow 0^+, \; \bar\delta=
\mu(\bar{\alpha}),\;\bar\alpha \in \sigma(TT^\ast),\;
c_1\,\bar\alpha\le\alpha_0, \;h(\bar{\alpha})<\norm{T}^2.
\end{equation}
Now, \allowdisplaybreaks
\begin{eqnarray}\label{eq:nueva1}\nonumber
  \norm{R_{\bar{\alpha}} Tx^\ast-x^\ast}^2 &=& \int_0^{\|T\|^{2\,+}} r_{\bar{\alpha}}^2(\lambda)\, d\norm{E_\lambda
  x^\ast}^2\\ \nonumber
   &\geq& \int_{h(\bar{\alpha})}^{\|T\|^{2\,+}} r_{\bar{\alpha}}^2(\lambda)\,d\norm{E_\lambda
   x^\ast}^2  ,\qquad \text{(since $0<h(\bar{\alpha})<\norm{T}^2$ )}   \\
   &\geq&\gamma^2 \rho^2(\bar{\alpha}) \int_{h(\bar{\alpha})}^{\|T\|^{2\,+}} s_\rho^{-2}(\lambda)\,d\norm{E_\lambda
   x^\ast}^2,
\end{eqnarray}
where the last inequality follows from the fact that
$(\rho,s_\rho)$ is an order-source pair for $\{g_\alpha\}$ (with
$\gamma$ the constant in (\ref{eq:4.50})\,). Since
$\eta(\bar{\delta})=\bar{\alpha}$ and $\bar{\delta}\to 0^+$ if and
only if $\bar{\alpha}\to 0^+$, (\ref{eq:nueva}) and
(\ref{eq:nueva1}) imply that
\begin{equation*}
 \int_{h(\bar{\alpha})}^{\|T\|^{2\,+}} s_\rho^{-2}(\lambda)\,d\norm{E_\lambda
   x^\ast}^2=o(1) \quad \textrm{for} \quad \bar{\alpha} \rightarrow 0^+
\end{equation*}
and therefore $\norm{s_\rho^{-1}(T^\ast T)x^\ast}=0$. Then
$x^\ast=0$, contradicting the fact that $x^\ast\neq 0$. Hence,
$\psi(x,\delta)=\rho(\Theta^{-1}(\delta))$ satisfies condition
\textit{S1} on $X^{s_\rho}.$

Also, since $\psi$ is trivially invariant over $X^{s_\rho}$, $\psi$ does not depend on $x$. Thus, it satisfies
condition \textit{S2}.

It only remains to be proved that $\psi$ satisfies condition
\textit{S3} on $X^{s_\rho}$. Suppose that is not the case. Then there must exist a set $M$,
$X^{s_\rho}\subsetneqq M\subset X\setminus \{0\}$ and $\tilde{\psi} \in
\mathcal{U}_{M}(\et)$ such that $\tilde{\psi}\mid
_{X^{s_\rho}}=\psi$ and $\tilde{\psi}$ satisfies \textit{S1} and
\textit{S2} on $M$. Let $x^\ast \in M\setminus X^{s_\rho}$. Since $\tilde{\psi} \in \mathcal{U}_{M}(\et)$ we have that
\begin{equation}\label{eq:des3} \et
\overset{\{x^\ast\}}{\preceq}\tilde{\psi}.
\end{equation}
Also, since $\tilde{\psi}$ is invariant over $M$ we have that
$\tilde{\psi}\overset{\{x^\ast\},X^{s_\rho}}{\preceq}{\tilde{\psi}}$,
and since $\tilde{\psi}$ coincides with $\psi$ on $X^{s_\rho}$, it
follows that
$\tilde{\psi}\overset{\{x^\ast\},X^{s_\rho}}{\preceq}{\psi}$.
This, together with (\ref{eq:des3}) implies that
$\et\overset{\{x^\ast\},X^{s_\rho}}{\preceq}{\psi}$ and therefore for every $x \in X^{s_\rho}$,
$\et(x^\ast,\delta)=O(\psi(x,\delta))$ as $\delta \rightarrow
0^+$, that is, $\et(x^\ast,\delta)=O(\rho(\Theta^{-1}(\delta)))$
as $\delta \rightarrow 0^+$. Lemma \ref{lem:o-O} then implies that
there exists an admissible \textit{a-priori} parameter choice rule
 $\tilde{\alpha}(\delta)$ such that
\begin{equation*}
    \underset{y^\delta \in \overline{B_\delta(Tx^\ast)}}{\sup}
\norm{R_{\tilde{\alpha}(\delta)}
y^\delta-x^\ast}=O(\rho(\Theta^{-1}(\delta))) \quad \textrm{as }
\delta \rightarrow 0^+
\end{equation*}
and therefore
\begin{equation*}
    \underset{y^\delta \in
    \overline{B_\delta(Tx^\ast)}}{\sup}\;\,\underset{\alpha\in (0,\alpha_0)}{\inf}
\norm{R_{\alpha} y^\delta-x^\ast}=O(\rho(\Theta^{-1}(\delta))) \quad
\textrm{as } \delta \rightarrow 0^+.
\end{equation*}
Hence, by virtue of Lemma \ref{lem:sup-inf}, $x^\ast \in \mathcal{R}(s_\rho(T^\ast T))$ and since $x^\ast \neq 0$, we have that $x^\ast \in \mathcal{R}(s_\rho(T^\ast T))\setminus \{0\}=X^{s_\rho}$ which contradicts our original assumption. This concludes the proof of the Theorem \ref{teo:sat-calopt}. \qed
\section{Examples}

Although the main results of this article are very theoretical in
nature, we provide below two examples of regularization methods
with optimal qualification which do possess saturation. In both
cases the saturation function and saturation set are found.

\begin{ej} The family of Tikhonov-Phillips regularization
operators $\{R_\alpha\}_{\alpha \in (0,\alpha_0)}$ where
$R_\alpha=g_\alpha(T^\ast T)T^\ast$ with
$g_\alpha(\lambda)=\frac{1}{\lambda + \alpha}$ has optimal
qualification $\rho(\alpha)=\alpha$
(\cite{ref:Herdman-Spies-Temperini-2009}). It can be easily checked
that $\rho$ is of local upper type 1,  $s_\rho(\lambda)=\lambda \in
\mathcal{S}$ and $\{g_\alpha\}_{\alpha \in (0,\alpha_0)}$ satisfies
all hypotheses of the Theorem \ref{teo:sat-calopt}. Therefore, the
function
$\psi(x,\delta)=\rho(\Theta^{-1}(\delta))=\delta^{\frac{2}{3}}$
defined for $x \in X^{s_\rho}\doteq \mathcal{R}(T^\ast
T)\setminus\{0\}$ and $\delta\in\left(0,\alpha_0^{\frac 3 2}\right)$
is saturation function of $\{R_\alpha\}_{\alpha \in (0,\alpha_0)}$
on $X^{s_\rho}$.
\end{ej}
\medskip

\begin{ej}
Given $k\in \mathbb{R}^+$, for $\alpha>0$ let
$$ u_\alpha^k(\lambda)\doteq \begin{cases}
\frac{e^{-\frac{\lambda}{\sqrt\alpha}}}{\lambda}, &\text{for } 0<\lambda<\alpha,\\
\frac{e^{-\sqrt\frac\lambda\alpha}}{\lambda}, &\text{for } \alpha\le\lambda<3\alpha,\\
\frac{e^{-\sqrt\frac\lambda\alpha}}{\lambda}+\frac{\alpha^k}{\lambda^{k+1}},
&\text{for } \lambda\ge 3\alpha,
\end{cases}
$$
and
$g_\alpha^k(\lambda)\doteq\frac{1}{\lambda}-\alpha^k\sqrt{\lambda}-u_\alpha^k(\lambda)$
for $\lambda>0$, and for $\lambda=0$ define
$g_\alpha^k(0)\doteq\underset{\lambda\to
0^+}{\lim}g_\alpha^k(\lambda)=\frac{1}{\sqrt{\alpha}}.$ The family
$\{g_\alpha\}_{\alpha\in (0,\alpha_0)}$ is a SRM
(\cite{ref:Herdman-Spies-Temperini-2011}).

Now, defining
$$ v_\alpha^k(\lambda)\doteq \begin{cases}
e^{-\frac{\lambda}{\sqrt\alpha}}, &\text{for } 0\leq\lambda<\alpha,\\
e^{-\sqrt\frac\lambda\alpha}, &\text{for } \alpha\le\lambda<3\alpha,\\
e^{-\sqrt\frac\lambda\alpha}+\left(\frac{\alpha}{\lambda}\right)^{k},
&\text{for } \lambda\ge 3\alpha,
\end{cases}$$
it follows that $r_\alpha^k(\lambda)\,=\,1-\lambda
g_\alpha^k(\lambda)\,=\,\alpha^k\lambda^{\frac 3
2}+v_\alpha^k(\lambda)$. If $\rho(\alpha)=\alpha^k$ then
$s_\rho(\lambda)=\lambda^k$, then  $s_\rho\in \mathcal{S}$ and
$s_\rho$ satisfies (\ref{eq:cond-calif}). Also, it can be show that
$\rho$ and $s_\rho$ verify (\ref{eq:4.50}) with $\gamma=1$  and
$h(\alpha)=3\alpha$. From Theorem \ref{teo:cond-calopt} it then
follows that $\rho(\alpha)=\alpha^k$ is optimal qualification of
$\{g_\alpha\}$.

On the other hand, for $k\geq 1$ and $\alpha>0$, the function
$g_\alpha^k(\lambda)$ is non-increasing. Thus, hypothesis
\textbf{b.v)} of Theorem \ref{teo:sat-calopt} holds and
$G_\alpha^k\doteq
\norm{g_\alpha^k(\cdot)}_\infty=g_\alpha^k(0)=\frac{1}{\sqrt{\alpha}}$,
which implies immediately that also hypothesis \textit{H4} is
verified. From now on we shall assume $k\ge 1$.

It can be easily checked that $\rho$ is of local upper type $k$
and $\{g_\alpha\}_{\alpha \in (0,\alpha_0)}$ satisfies all
hypotheses of the Theorem \ref{teo:sat-calopt}. Therefore, the
function
$\psi(x,\delta)=\rho(\Theta^{-1}(\delta))=\delta^{\frac{2k}{2k+1}}$
defined for $x \in X^{s_\rho}\doteq \mathcal{R}((T^\ast
T)^k)\setminus\{0\}$ and $\delta\in\left(0,\alpha_0^{k+\frac 1
2}\right)$ is saturation function of $\{R_\alpha\}_{\alpha \in
(0,\alpha_0)}$ on $X^{s_\rho}$.
\end{ej}

\section{Conclusions}

In this article families of real functions $\{g_\alpha\}$ defining
a spectral regularization methods with optimal qualification were
considered. Sufficient conditions on the family and on the optimal
qualification guaranteeing the existence of saturation were found.
Appropriate characterizations of both the saturation function and
the saturation set were given and two examples were provided.




\end{document}